\catcode`\@=11

\magnification=1200
\baselineskip=14pt

\pretolerance=500    \tolerance=1000 \brokenpenalty=5000

\catcode`\;=\active
\def;{\relax\ifhmode\ifdim\lastskip>\z@
\unskip\fi\kern.2em\fi\string;}

\overfullrule=0mm

\catcode`\!=\active
\def!{\relax\ifhmode\ifdim\lastskip>\z@
\unskip\fi\kern.2em\fi\string!}

\catcode`\?=\active
\def?{\relax\ifhmode\ifdim\lastskip>\z@
\unskip\fi\kern.2em\fi\string?}

\frenchspacing

\newif\ifpagetitre            \pagetitretrue
\newtoks\hautpagetitre        \hautpagetitre={\hfill \date }
\newtoks\baspagetitre         \baspagetitre={1}

\newtoks\auteurcourant        \auteurcourant={   }
\newtoks\titrecourant
\titrecourant={  }

\newtoks\hautpagegauche       \newtoks\hautpagedroite
\hautpagegauche={\hfill\sevenrm\the\auteurcourant\hfill}
\hautpagedroite={\hfill\sevenrm\the\titrecourant\hfill}

\newtoks\baspagegauche       \baspagegauche={\hfill\rm\folio\hfill}

\newtoks\baspagedroite       \baspagedroite={\hfill\rm\folio\hfill}

\headline={
\ifpagetitre\the\hautpagetitre
\global\pagetitrefalse
\else\ifodd\pageno\the\hautpagedroite
\else\the\hautpagegauche\fi\fi}

\footline={\ifpagetitre\the\baspagetitre
\global\pagetitrefalse
\else\ifodd\pageno\the\baspagedroite
\else\the\baspagegauche\fi\fi}

\def\date{\ {\the\day}\
\ifcase\month\or Janvier\or F\'evrier\or Mars\or Avril
\or Mai \or Juin\or Juillet\or Ao\^ut\or Septembre
\or Octobre\or Novembre\or D\'ecembre\fi\
{\the\year}}

\def\up#1{\raise 1ex\hbox{\sevenrm#1}}

\def\cqfd{\unskip\kern 6pt\penalty 500
\raise -2pt\hbox{\vrule\vbox to 10pt{\hrule width 4pt
\vfill\hrule}\vrule}\par\medskip}

\def\section#1{\vskip 7mm plus 20mm minus 1.5mm\penalty-50
\vskip 0mm plus -20mm minus 1.5mm\penalty-50
{\bf\noindent#1}\nobreak\smallskip}

\def\subsection#1{\medskip{\bf#1}\nobreak\smallskip}

\def\displaylinesno #1{\dspl@y\halign{
\hbox to\displaywidth{$\@lign\hfil\displaystyle##\hfil$}&
\llap{$##$}\crcr#1\crcr}}

\def\ldisplaylinesno #1{\dspl@y\halign{
\hbox to\displaywidth{$\@lign\hfil\displaystyle##\hfil$}&
\kern-\displaywidth\rlap{$##$}
\tabskip\displaywidth\crcr#1\crcr}}

\def\hfl#1#2{\smash{\mathop{\hbox to 12 mm{\rightarrowfill}}
\limits^{\scriptstyle#1}_{\scriptstyle#2}}}

\catcode`\@=12

\def\bR{{\bf R}}
\def\bZ{{\bf Z}}

\def\cU{{\cal U}}
\def\cV{{\cal V}}

\def\ux{{\underline{ x}}}

\def\uy{{\underline{ y}}}

\def\uz{{ \underline{z}}}

\def\bZ{{\bf Z}}

\def\utheta{{\underline{\theta}}}

\def\proof{\bigskip\noindent{\it Proof.}\ }

\def\and{\quad\hbox{and}\quad}


\def\M{\mathop{\rm M\kern 1pt}\nolimits}
\def\h{\mathop{\rm h\kern 1pt}\nolimits}

\def\romain#1{\uppercase\expandafter{\romannumeral #1}}

\def\card{\mathop{\rm Card\kern 1.3 pt}\nolimits}
\def\deg{\mathop{\rm deg\kern 1pt}\nolimits}
\def\det{\mathop{\rm det\kern 1pt}\nolimits}

\def\h{\mathop{\rm h\kern 1pt}\nolimits} \long\def\forget#1\endforget{}

\def\og{\leavevmode\raise.3ex\hbox{$\scriptscriptstyle
\langle\!\langle\,$}}
\def\fg{\leavevmode\raise.3ex\hbox{$\scriptscriptstyle
\!\rangle\!\rangle\,\,$}}


\def\ArRo{1}
\def\BuGl{2}
\def\BuLiv{3}
\def\BuCh{4}
\def\BuLa{5}
\def\Cas{6}
\def\CasB{7}
\def\DaSc{8}
\def\DSa{9}
\def\FalA{10}
\def\FalB{11}
\def\Gro{12}
\def\Jar{13}
\def\JarB{14}
\def\Khi{15}
\def\LagA{16}
\def\LagB{17}
\def\Lau{18}
\def\RoyA{19}
\def\RoyB{20}
\def\ScTr{21}
\def\Schm{22}
\def\Wa{23}


\centerline{}

\vskip 4mm

\centerline{
\bf Exponents of inhomogeneous Diophantine Approximation}

\vskip 8mm
\centerline{Yann B{\sevenrm UGEAUD} \footnote{}{\rm
2000 {\it Mathematics Subject Classification : } 11J20, 11J13, 11J82.}
\&
Michel L{\sevenrm AURENT} }

\vskip 9mm

\noindent {\bf Abstract --} In Diophantine Approximation,
inhomogeneous problems are linked with homogeneous ones by means of
the so-called Transference Theorems. We revisit this classical topic by
introducing new exponents of Diophantine approximation. We prove that 
the
exponent of approximation to a generic point in $\bR^n$ by a system of 
$n$ linear forms
is equal to the inverse of the uniform homogeneous exponent associated
to the system of dual linear forms.

\vskip15mm

\section
{1. Introduction and results.}

It is a well known fact that inhomogeneous problems
in Diophantine Approximation are connected to homogeneous
ones by means of the so-called Transference Theorems. We revisit
this classical topic, refering mainly to the book of Cassels [\Cas],
in the light of
recent results on  uniform exponents of Diophantine approximation,
which have been introduced in a restricted setting  in [\BuLa].
These uniform exponents, indicated by a `hat', enable us here to control
quantitatively the accuracy of the
inhomogeneous approximation at a generic point in the ambient space.

Let us begin with some notations and definitions.
If  $\utheta$ is a  (column) vector  in $\bR^n$, we denote by $\vert
\utheta\vert$
the  maximum of the   absolute values of its coordinates  and by
$$
\Vert \utheta\Vert = \min_{\ux \in \bZ^n} \vert \utheta -\ux\vert
$$
the maximum of the distances of its  coordinates to the  rational
integers.

Let $n$ and $m$ be two positive integers and let  $A$ be a real matrix
with $n$ rows and $m$ columns.
     For any $n$-tuple $\utheta$
of real numbers, we denote by  $w(A,\utheta)$
the supremum of the real numbers $w$ for which, {\it for arbitrarily
large real numbers $X$}, the inequalities
$$
\Vert A\ux + \utheta\Vert \le X^{-w} \and \vert \ux \vert \le X
\leqno{(1)}
$$
have a solution  $\ux$ in $\bZ^m$. Following the conventions of
[\BuLa], we
denote by  $\hat{w}(A,\utheta)$ the supremum of the real numbers $w$
for which, {\it for all sufficiently large positive real numbers  $X$},
the inequalities (1) have an integer solution $\ux$ in $\bZ^m$.
The lower bounds
$$
w(A,\utheta) \ge \hat{w}(A,\utheta) \ge 0
$$
are then obvious. We define furthermore
    two homogeneous  exponents
     $w(A)$ and $\hat{w}(A)$ as in  (1) with $\utheta
= {}^t (0,\dots , 0)$,
    requiring moreover  that the integer solution $\ux$  should be
non-zero.
Then the Dirichlet box principle implies that
$$
w(A) \ge  \hat{w}(A) \ge {m\over n}.
$$
In this respect, Jarnik [\JarB] has given non trivial lower bounds
of $w(A)$ in terms of $\hat{w}(A)$.
Furthermore, it follows from the Borel--Cantelli Lemma    that the
equalities
$$
w(A) =  \hat{w}(A) = {m\over n}
$$
hold for almost all matrices $A$ in $M_{n, m}(\bR)$
with respect to the Lebesgue measure on $\bR^{mn}$
(see also
Groshev [\Gro]; his result is sharper, actually).

The transposed  matrix of any  matrix
$A$ is denoted by ${}^t A$. We can now state our main result.

\proclaim
Theorem. For any $n$-tuple $\utheta$ of real numbers, we have the lower
bounds
$$
w(A,\utheta) \ge {1\over \hat{w}({}^tA)} \and
\hat{w}(A,\utheta) \ge {1\over {w}({}^tA)} ,
\leqno({2})
$$
with equality  in (2) for almost all $\utheta$ with respect to the
    Lebesgue measure on $\bR^n$.

Let us first examine the simplest case of a $1\times 1$ matrix
$A=(\xi)$, where
$\xi$ is an irrational real number.  The Dirichlet box principle
asserts that
for any real number $Q\ge 1$, there exists an integer  $q$ such that
$$
1 \le \vert q\vert \le Q \and \Vert q\xi \Vert \le Q^{-1} .
$$
This statement is best possible in the sense that
the exponent $-1$ in the upper bound $Q^{-1}$ cannot be replaced
by a smaller value (see e.g. [\Khi],
[\DSa] or [\Schm], page 26). It then follows that
our exponent  $\hat{w}({}^tA)$ is equal to  $1$, as well as the
generic inhomogeneous exponent $\displaystyle\inf_{\theta\in
\bR}w(A,\theta)$.
In fact, a more precise result holds in this particular case.
Namely, Minkowski has proved  that for any real number
$\theta$, the system of inequalities
$$
\vert q \vert \le Q \and \Vert q\xi +\theta\Vert \le {1\over 4} Q^{-1}
$$
has an integer solution
    $q$ for infinitely many integers $Q$. Moreover, our Theorem
    shows that $-1$  is the best
possible exponent, regardless of the irrational number $\xi$.
Besides, Cassels has observed that there does not exist any
inhomogeneous analogue of
the Dirichlet box principle, even if we weaken the property of
approximation.
In Theorem III of Chapter 3 from [\Cas], he constructed
a Liouville number $\xi$ and a real number
$\theta $ such that, for any  $\epsilon >0$,  we have the lower bound
$$
     \min_{\vert q\vert \le Q}\Vert q\xi +\theta\Vert \ge Q^{-\epsilon}
$$
for infinitely many integers $Q$. The result of Cassels follows
immediately
from our
    Theorem, since $w((\xi))=+\infty$ for any
    Liouville number $\xi$. The
     uniform exponent of inhomogeneous approximation
$\hat{w}((\xi),\theta)$ therefore
vanishes for almost all $\theta$.

\medskip

The cases $(m, n) = (1, 2)$ and $(m, n) = (2, 1)$ have been
studied by Khintchine [\Khi]. Notice that the metrical statement of our
Theorem is of same spirit as Satz V b from [\Khi].

\medskip

We could refine the first statement of the Theorem
(and of Proposition 1 below) by taking into account
whether or not there exists a positive constant $c$ such that,
{\it for arbitrarily large real numbers $X$}, the inequalities
$$
\Vert A\ux + \utheta\Vert \le c \, X^{-w(A, \utheta)}
\and \vert \ux \vert \le X
$$
have a solution  $\ux$ in $\bZ^m$. We take into consideration
this remark in the statement of the Corollary below.

\medskip

The next result should be compared with Theorem X of
Chapter 5 from [\Cas].

\proclaim
Proposition 1.
Let  $A$ be a  matrix in $M_{n,m}(\bR)$.
For any exponent $w> 1/ \hat{w}({}^tA)$,
there exists a real $n$-tuple $\utheta$
such that the lower bound
$$
\Vert A\ux + \utheta\Vert \ge \vert\ux\vert^{-w}
$$
holds for any integer $m$-tuple $\ux$ whose norm  $\vert\ux\vert$
is sufficiently
large. Moreover, there exists some real  $n$-tuple $\utheta$ such that
$$
\Vert A\ux + \utheta\Vert \ge
{1 \over 72 n^2 (8m)^{m/n}} \, \vert\ux\vert^{-m/n}
$$
holds for all non-zero integer $m$-tuples $\ux$.

Cassels [\Cas], page 85, has proved the second assertion of Proposition
1,
without however computing the
constant occurring in the right-handside of the lower bound.
Moreover, our first  assertion improves
    Theorem X of  [\Cas] whenever $\hat{w}({}^tA) > n/m$. See
Section 5 for examples of such matrices $A$.

If  the subgroup $G= {}^tA\bZ^n + \bZ^m$ of $\bR^m$ generated
by the $n$ rows of $A$ together with $\bZ^m$
has maximal rank $m+n$,  then Kronecker's Theorem
asserts that the dual subgroup $\Gamma = A \bZ^m +\bZ^n$ of $\bR^n$
generated by the
$m$ columns of  $A$ and by $\bZ^n$  is dense in
$\bR^n$. With this respect, our Theorem may be viewed as a  measure of
the density of $\Gamma$.
In the case where the rank of  $G$ is
$< m+n$, we clearly have
$$
\hat{w}({}^t A) = w({}^t A) =+\infty \and
\hat{w}(A,\utheta)={w}(A,\utheta)=0
$$
for any $n$-tuples  $\utheta$ located outside a discrete family of
parallel hyperplanes  in $\bR^n$. The
assertion of the  Theorem is then obvious.
In the sequel of the paper, we shall therefore assume that
the rank over $\bZ$ of the group $G$ is equal to $m+n$.
Notice however that the exponent $\hat{w}({}^t A)$
may be  infinite, even when $G$ has  rank $m+n$, as
proved by Khintchine [\Khi] in the case $(m, n) = (1, 2)$
and by Jarn\'\i k [\Jar] in the general case $(m, n)$
with $n \ge 2$ (when $n=1$ and $m$ is arbitrary,
the exponent $\hat{w}({}^t A)$ can be as large as 1 [\Jar], but
not larger [\Khi]). See also Theorem XIV (page 94)
of [\Cas] for the construction of such a  matrix
$A$ and the following Theorem XV concerning the density of the
associated group $\Gamma$.

Let us illustrate our Theorem by the example of the row (resp. column)
matrices
$$
A = (\xi, \dots , \xi^n) , \quad  \hbox{\rm resp.} \quad A =
{}^t(\xi, \dots , \xi^n) ,
$$
made up with the successive powers of a transcendental real number
$\xi$.
Then, the corresponding exponents $\hat{w}(A)$ are  uniformly bounded
in terms of $n$ (see [\BuLa] for references).
Roy [\RoyA] determined these exponents for $n=2$ when
$\xi$ is a Fibonacci continued fraction, that is, when we have
$$
\xi = [0; a, b, a, a, b, a, b, \ldots],
$$
where the sequence of partial quotients of $\xi$ is given by the
fixed point of the Fibonacci substitution $a \to ab$, $b \to a$. Here,
$a$ and $b$ denote distinct positive integers.

Combining these results with our
Theorem, we obtain the following statement.

\bigskip

\noindent
\bf
Corollary. \sl Let $n$ be a positive integer and let $\xi$ be a
real transcendental number.

(i) There exists a positive constant $c$ such that, for any real number
$\theta$, there exist infinitely many
polynomials $P(X)$ with integer coefficients, degree
at most $n$, and
$$
\vert P(\xi)+\theta\vert \le c H(P)^{-\lceil n/2 \rceil}.
$$

(ii)  There exists a positive constant $c$ such that, for
any  real $n$-tuple  $\utheta = (\theta_1,\dots , \theta_n)$,
there exist infinitely many integers $q$ with
$$
\max_{1\le j\le n}\Vert q\xi^j +\theta_j\Vert \le
c \, \vert q\vert ^{-1/(2n-1)}.
$$

(iii) When  $n=2$, the assertions (i) and  (ii)
remain  valid with the exponents
$$
(1+\sqrt{5})/2 \simeq 1.618... \and
(3-\sqrt{5})/2\simeq 0.3819...,
$$
respectively. If moreover  $\xi$ is a  Fibonacci continued fraction
and if
$$
w > (1+\sqrt{5})/2 \and  \lambda > (3-\sqrt{5})/2
$$
then, for almost all real numbers
$\theta$ and for almost all  pairs of real numbers
$(\theta_1, \theta_2)$, we have the respective lower bounds
$$
\vert P(\xi)+\theta\vert \ge H(P)^{-w}
\and
\max \{\Vert q\xi +\theta_1 \Vert, \Vert q\xi^2 +\theta_2 \Vert\}
\ge \vert q\vert^{- \lambda},
$$
for any polynomial $P(X)$ with integer coefficients and degree $\le 2$
and with sufficiently large height, and for any integer $q$ with
sufficiently large absolute value.

\rm
\bigskip

Thus, for a Fibonacci continued fraction  $\xi$, the critical exponents
in degree $n=2$ for the  problems of inhomogeneous approximation
    (i) and (ii) are respectively $w=1.618...$
and  $\lambda =  0.3819...$,
instead of the exponents $2$ and $1/2$, which occur in the generic 
situation.
Notice that $\hat{w}(A)$ and $\hat{w}({}^t A)$ have also been
determined for $A = (\xi, \xi^2)$ when $\xi$ is a Sturmian continued
fraction, see [\BuLa]. This provides further examples of matrices $A$
with $\hat{w}(A)$ and $\hat{w}({}^t A)$ greater than in the
generic case.

\medskip

To conclude this introduction, let us recall the inequalities
relying the exponents $w(A)$ and $w({}^t A)$, which follow from the
    Khintchine transference principle (cf. for example [\Schm],
Theorem 5C, page 99)\thinspace ;
with the preceding notations, we have
$$
w(A) \ge {m \, w({}^t A) + m - 1 \over (n - 1) w({}^t A) + n}.
$$
Furthermore, a careful reading of the proof shows that the uniform
exponents
$\hat{w} (A)$ and  $\hat{w} ({}^t A)$ are linked by the same relation
$$
\hat{w} (A) \ge {m \, \hat{w} ({}^t A) + m - 1
\over (n-1) \hat{w} ({}^t A) + n}.
$$

\medskip

Our article is organized as follows.
Section 2 is devoted to the definition and to the properties of the
sequence of best approximations.
A crucial fact  is that it increases at least geometrically.
Some transference lemma is stated and proved in  Section 3. It is used
in the
next Section, where we establish the Theorem and Proposition 1.
    The Corollary  is then  discussed in Section 5. Finally,
    questions of  Hausdorff dimensions, which arise naturally from the
Theorem, are briefly treated in Section 6.

\section
{2. Best approximations.}

Following the notations of  [\Cas], we denote by
$$
M_j(\uy) = \sum_{i=1}^{n}\alpha_{i,j}y_i , \quad \uy ={}^t (y_1, \dots
, y_n),
\quad (1 \le j \le m)
$$
the linear forms  determined by the  columns of the  matrix
$A= (\alpha_{i,j})$ and we set
$$
M(\uy) = \max_{1\le j\le m}\Vert M_j(\uy)\Vert = \Vert {}^tA\uy\Vert.
$$
Observe that the quantity $M(\uy)$ is  positive for all non-zero integer
    $n$-tuples $\uy$,
since we have assumed that the rank over
$\bZ$ of the group $G$ is equal to  $m+n$.
Then, we can build inductively a sequence of integer vectors
$$
\uy_i = {}^t(y_{i,1}, \dots , y_{i,n}) ,\qquad (i\ge 1),
$$
called {\it a sequence of  best approximations}
\footnote{(*)}{ According to [\LagB], a best approximation should be
     a vector belonging to  $\bZ^{n+m}$, by analogy with the usual
continued fraction process.
We forget here the last $m$
coordinates which are unsignificant  for our purpose. }
related to the linear  forms
$M_1, \dots , M_m$ and to the supremum  norm,  which satisfies the
following properties.  Set
$$
\vert \uy_{i }\vert = Y_i  \and  M_i =  M(\uy_i).
$$
Then, we have
$$
1 = Y_1 < Y_2<\cdots  \and   M_1 > M_2 > \cdots \enspace ,
$$
and $M(\uy)\ge M_i$ for all non-zero integer vectors   $\uy$ of norm
$\vert \uy
\vert < Y_{i+1} $.
We start the construction
with a  smallest {\it  minimal point } $\uy_1$ in the sense of  [\DaSc],
verifying
$Y_1 = \vert \uy_{1}\vert = 1$
     and   $M(\uy)\ge M(\uy_1) =M_1$
for any  integer point  $\uy \in \bZ^n$ with norm $\vert \uy
\vert = 1$.
Suppose that   $\uy_1, \dots , \uy_i$ have already been constructed
in such a way that $M(\uy)\ge M_i$ for all non-zero integer point $\uy$
of norm
$\vert\uy\vert \le Y_i$. Let  $Y$
be the smallest positive integer  $>Y_i$ for which there exists an
integer point $\uz$
verifying  $\vert \uz \vert =Y$ and  $ M(\uz) < M_i$.
The integer  $Y$ does exist by the Dirichlet box principle since
$M_i>0$.
Among those points $\uz$, we select an element $\uy$ for which
$M(\uz)$ is minimal.
   We then set
$$
\uy_{i+1} = \uy, \qquad
Y_{i+1} = Y \and M_{i+1} =  M(\uy).
$$
The sequence  $(\uy_i)_{i\ge 1}$ thus obtained  satisfies clearly the
desired properties.

Let  $w$ be a  real number $ <\hat{w}({}^t A)$, so that
the system of inequalities
$$
M(\uy) \le Y^{-w} \and \vert\uy\vert \le Y
$$
has a non-zero integer solution $\uy$ for any sufficiently large $Y$.
Choosing $Y<Y_{i+1}$  arbitrarily close to  $Y_{i+1}$, we obtain the
upper bound
$$
M_i \le Y_{i+1}^{-w} \leqno{(3)}
$$
for any sufficiently large index $i$,  using the  characteristic
property
of the best approximations.

Suppose now that $w<w({}^t A)$. Then, there exist
infinitely many indices $i$ for which
$$
M_i \le Y_{i}^{-w} . \leqno{(4)}
$$
The  indices $i$ satisfying  (4) are obtained by  inserting the  norm
$$
Y_i\le \vert\uy\vert <Y_{i+1}
$$
of the integer solutions  $\uy$ of the inequation $M(\uy)\le
\vert\uy\vert^{-w}$
in the sequence  $(Y_k)_{k\ge 1}$.

Observe furthermore that the Dirichlet box principle
(cf. [\Cas], Theorem VI, page 13) ensures that the system of
inequations
$$
M(\uy) \le Y^{-n/m} \and \vert\uy\vert \le Y
$$
has a non-zero integer solution  $\uy$ for any  $Y \ge 1 $.
Arguing as above, we obtain the upper bound
$$
M_i \le Y_{i+1}^{-n/m} \leqno{(5)}
$$
for all $i \ge 1$.

     \proclaim
Lemma 1.  There exists a positive constant
    $c$ such that
$$
Y_i \ge c \, 2^{i/(3^{m+n}-1)}
$$
for all $i\ge 1$.

\proof  Lagarias [\LagB ] has established in a quite  general framework
that a sequence of best
    approximations increases at least
geometrically. We take again the argumentation used in the proof of
Theorem 2.2 from [\LagA].
Let us  consider the
$3^{m+n} +1$ consecutive vectors
$$
\uy_i, \uy_{i+1}, \dots , \uy_{i+3^{m+n}}.
$$
By the usual box principle, there exist two  indices $r$ and $s$, with
$0\le r < s \le 3^{m+n}$, such that
$$
\uy_{i+r,j} \equiv \uy_{i+s,j} \quad ({\rm mod}\,3) \quad \hbox{\rm for
all } j=1,\dots ,n,
$$
and
$$
\langle M_k(\uy_{i+r})\rangle \equiv \langle M_k(\uy_{i+s})\rangle
\quad ({\rm mod}\, 3) \quad \hbox{\rm for  all } k=1, \dots , m,
$$
where the notation $\langle x\rangle$ stands for the closest integer
to the real number $x$.
   Setting
$$
\uz ={\uy_{i+s}-\uy_{i+r}\over 3},
$$
   we have
$$
\vert \uz \vert \le {Y_{i+s}+Y_{i+r}\over 3}
\and
M(\uz) \le {M_{i+r}+M_{i+s}\over 3} < M_{i+r}.
$$
Since $\uz$ is a non-zero integer vector, we get
$$
Y_{i+r+1} \le {Y_{i+s} + Y_{i+r}\over 3}
$$
and
$$
Y_{i+3^{m+n}}\ge Y_{i+s}\ge 3 Y_{i+r+1} -Y_{i+r} \ge 2 Y_{i+r +1}\ge 2
Y_{i+1}
$$
for any $i\ge 1$. The expected lower bound then follows  by induction
on $i$. \cqfd

\bigskip

\proclaim
Lemma 2. For almost all real $n$-tuples  $\utheta = {}^t(\theta_1, \dots
,\theta_n)$, we have the lower bound
$$
\Vert y_{i,1}\theta_1 + \cdots + y_{i,n}\theta_n\Vert  \ge
Y_i^{-\delta},
$$
     for any  $\delta >0$,
and any index  $i$ which is sufficiently large in terms of $\delta$
and of  $\theta_1 , \dots ,\theta_n$.

\proof
We can assume without restriction that the numbers  $\theta_j$ are
located in the interval
$[0,1]$ and that $\delta $
is given. For a fixed   $i$, the reverse inequality
$$
\Vert y_{i,1}\theta_1 + \cdots + y_{i,n}\theta_n\Vert  <
Y_i^{-\delta} \leqno{(6)}
$$
defines in the  hypercube $[0,1]^n$ a subset of Euclidean  volume
    $\le 2\sqrt{n}Y_i^{-\delta}$. By Lemma 1,
the series $\sum_{i\ge 1} Y_i^{-\delta}$ converges.
It then follows from the
     Borel--Cantelli Lemma that the set of  $\utheta$ satisfying $(6)$ 
for
    infinitely many  indices $i$
has Lebesgue measure zero. \cqfd

\section
{3. A transference lemma.}
Let us consider now the linear forms
$$
L_i(\ux) = \sum_{j=1}^{m}\alpha_{i,j}x_j , \quad \ux ={}^t(x_1, \dots ,
x_m),
\quad (1 \le i \le n),
$$
     determined by the rows of the matrix $A$.
The following result relies the  problem of the inhomogeneous
simultaneous approximation
by the linear forms $L_i$ to the problem of
homogeneous simultaneous approximation by the linear forms $M_j$.

\proclaim
Lemma 3. Set  $\kappa = 2^{1-m-n}((m+n)!)^2$. Let  $X$ and $Y$ be
two positive real numbers.
Suppose that we have the lower bound
$$
M(\uy) \ge \kappa X^{-1}
$$
     for any non-zero integer $n$-tuple $\uy$ of
    norm $\vert\uy\vert \le Y$. Then, for all real  $n$-tuples
$(\theta_1, \dots , \theta_n)$, there exists an integer
$m$-tuple $\ux$ with norm $\vert\ux\vert \le X $ such that
$$
     \max_{1\le i\le n}\Vert L_i(\ux) + \theta_i \Vert \le \kappa Y^{-1}.
$$

\proof
This is the first assertion of Lemma 4.1 from [\Wa]. For the
convenience of the reader,
we reproduce the proof.
Let $X$ and $Z$ be two positive real numbers. Part B of
    Theorem
XVII in chapter V of [\Cas] asserts that the    system of  inequations
$$
\max_{1\le i\le n}\Vert L_i(\ux) + \theta_i \Vert \le Z \and
\vert\ux\vert \le X
$$
has an integer solution $\ux\in \bZ^m$, provided that the upper bound
$$
\Vert y_1 \theta_1  + \cdots + y_n\theta_n \Vert \le
\kappa^{-1}\max\{XM(\uy),Z\vert\uy\vert\}
\leqno{(7)}
$$
holds for all integer $n$-tuples $\uy$. Let us apply this result with
$Z= \kappa Y^{-1}$.
The condition (7) is satisfied when $\vert\uy\vert \ge Y$, since then
$\kappa^{-1}Z\vert\uy\vert$ is $\ge 1$ while the left hand side of (7)
is $\le 1/2$. If
$\uy$ is non-zero  and $\vert\uy\vert \le Y$, our assumption ensures
that
    $M(\uy)\ge \kappa X^{-1}$,
and the right hand side of  (7) is  $\ge 1$ in this case too. Finally,
(7) obviously holds  for $\uy=0$.\cqfd

\medskip
Up to the value of the numerical constant  $\kappa$, the above
mentioned Theorem XVII of [\Cas]
is a consequence of the more general Theorem VI in Chapter XI of
[\CasB], when applied to the distance function
$$
F(x_1, \dots , x_{m+n}) = X^{-1}\biggl( \, \sum_{j=1}^m\vert x_j\vert 
\, \biggr)
+Z^{-1} \biggl( \,  \sum_{i=1}^n \vert L_i(x_1,\dots , x_m)
+x_{m+i}\vert\, \biggr)
$$
in $\bR^{m+n}$.
Notice that this last result provides also an explicit construction
of the approximating point $\ux$ in terms of successive minima and of
duality.

\section
{4. Proof of the Theorem and of Proposition 1.}

First, we prove that the lower bounds
$$
{w}(A,\utheta) \ge {1\over \hat{w}({}^tA)} \and
\hat{w}(A,\utheta) \ge {1\over {w}({}^tA)}
\leqno{(8)}
$$
hold for all real $n$-tuples $\utheta ={}^t(\theta_1,
\dots , \theta_n)$.

Let $w > \hat{w}({}^tA)$ be a real number. By definition of the exponent
$\hat{w}({}^tA)$, there exists
a real number $Y$,  which may be chosen arbitrarily large, such that
$$
M(\uy) \ge Y^{-w} \leqno {(9)}
$$
for any non-zero integer $n$-tuple $\uy$ of norm $\vert \uy\vert
\le Y$. We use
    Lemma 3 with $X = \kappa Y^w$, where
$\kappa = 2^{1-m-n}((m+n)!)^2$. Consequently, there exists an integer
$m$-tuple
    $\ux$ of norm $\vert\ux\vert \le   X$ such that
$$
     \max_{1\le i\le n}\Vert L_i(\ux) + \theta_i \Vert \le \kappa Y^{-1}
=\kappa^{(1+1/w)}X^{-1/w}
\le \kappa^{(1+1/w)}\vert\ux\vert^{-1/w}.
$$
We deduce that  $w(A,\utheta) \ge 1/w$. The first assertion
of  (8) then follows by letting $w$ tend to  $\hat{w}({}^tA)$.

The second lower bound of  (8) is established along the same lines,
observing that for
$w > {w}({}^tA)$ and any sufficiently large real number $Y$,
inequality (9) is satisfied for any
non-zero integer $n$-tuple $\uy$ of norm $\vert \uy\vert\le Y$.

We shall now prove that the inverse upper bounds
$$
{w}(A,\utheta) \le {1\over \hat{w}({}^tA)} \and
\hat{w}(A,\utheta) \le {1\over {w}({}^tA)}
\leqno{(10)}
$$
hold for almost all real  $n$-tuples $\utheta=
{}^t(\theta_1, \dots , \theta_n)$.

The duality formula  ${}^t \uy A \ux = {}^t \ux {}^t A \uy$  written
in the form
$$
y_1\theta_1 + \cdots + y_n\theta_n = \sum_{i=1}^n y_i\left(L_i(x_1,
\dots, x_m)+\theta_i\right)
-\sum_{j=1}^mx_jM_j(y_1,\dots , y_n)
$$
    implies  the upper bound
$$
\Vert y_1\theta_1 + \dots + y_n\theta_n \Vert
\le n\vert \uy \vert \max_{1\le i\le n}\Vert L_i(\ux) + \theta_i \Vert
+ m \vert \ux \vert M(\uy) \leqno{(11)}
$$
for all integer vectors  $\ux={}^t(x_1, \dots, x_m)$ and $\uy=
{}^t(y_1, \dots, y_n)$.

Let
$$
\uy_i = {}^t(y_{i,1}, \cdots , y_{i,n}) \and Y_i = \vert\uy_i\vert
\quad (i\ge 1)
$$
be a sequence of best approximations relative to the matrix
${}^tA$. Suppose that
for all  $\delta >0$ we have the lower bound
$$
\Vert y_{i,1}\theta_1 + \cdots + y_{i,n}\theta_n\Vert \ge Y_i^{-\delta}
\leqno{(12)}
$$
for any index  $i$  large enough. By Lemma 2,
the estimation (12)
holds for almost all real  $n$-tuple $\utheta$. Let us fix now
two real numbers  $\delta$ and  $w$ such that
$$
0 < \delta < w < \hat{w}({}^tA).
$$
Let  $\ux$ be an integer  $m$-tuple  with sufficiently large norm
$\vert\ux\vert$,
    and let $k$ be the index defined by the inequalities
$$
Y_k \le (2m\vert\ux\vert)^{1/(w - \delta)} < Y_{k+1},
$$
so that
$$
Y_{k+1}^w > (2m \vert \ux\vert)^{w/ (w -\delta)} \ge 2m\vert\ux\vert
Y_k^\delta.
$$
Combining now the estimations (3),
     (11) with $\uy=\uy_k$ and (12) for $i=k$, we obtain
$$
\eqalign{
Y_k^{-\delta}  \le & \enspace
n \, \vert \uy_k \vert \, \max_{1\le i\le n}\Vert L_i(\ux) + \theta_i
\Vert
+ m \, \vert \ux \vert \, M(\uy_k)
\cr
\le & \enspace n \, Y_k \, \max_{1\le i\le n} \Vert L_i(\ux) + \theta_i
\Vert
+ m \, \vert \ux \vert \, Y_{k+1}^{-w}
\cr
\le & \enspace n \, Y_k \, \max_{1\le i\le n} \Vert L_i(\ux) + \theta_i
\Vert
+ {Y_k^{-\delta} \over 2},
\cr}
$$
from which follows the  lower bound
$$
\eqalign{
\Vert A\ux +\utheta\Vert & =
     \max_{1\le i\le n}\Vert L_i(\ux) + \theta_i \Vert
\ge {1\over 2n} Y_k^{-1-\delta} \cr
& \ge  (2m)^{-(\delta +1)/ (w -\delta)}(2n)^{-1}
\vert \ux\vert^{-(\delta +1)/ (w -\delta)}. \cr}
$$
We deduce that
$$
w(A,\utheta) \le {\delta +1\over w-\delta}.
$$
Choosing $\delta$ and $w$ arbitrarily close to  $0$ and to
$\hat{w}({}^tA)$
respectively, we obtain the first  upper bound of (10).

In order to prove the second upper bound of (10), we take again the
preceding argumentation
using now the estimation (4) instead of (3).
Let us fix two real numbers  $\delta$ and $w$ satisfying
$$
0 < \delta < w < {w}({}^tA).
$$
Let  $k$ be an integer  $\ge 1$ such that $M_k \le Y_k^{-w}$.
Since $w < w({}^tA)$,
there exist infinitely many such integers $k$. Let $\ux$ be an
integer $m$-tuple with norm
$\vert\ux\vert \le X_k := Y_k^{w-\delta}/ (2m)$.  Combining
(4), (11) and  (12), we obtain
$$
\eqalign{
Y_k^{-\delta}  \le & \enspace
n \, \vert \uy_k \vert \, \max_{1\le i\le n}\Vert L_i(\ux) + \theta_i
\Vert
+ m \, \vert \ux \vert \, M(\uy_k)
\cr
\le & \enspace n \, Y_k \, \max_{1\le i\le n} \Vert L_i(\ux) +
\theta_i \Vert + m \, X_k \, Y_k^{-w},
\cr}
$$
from which we deduce that
$$
     \max_{1\le i\le n}\Vert L_i(\ux) + \theta_i \Vert
\ge {1\over 2n} Y_k^{-1-\delta} =  (2m)^{-(\delta +1)/ (w
-\delta)}(2n)^{-1}
X_k^{-(\delta +1)/ (w -\delta)}.
$$
Recall now that the above lower bound holds for any integer point
with   norm $\le X_k$
and for infinitely many integers $k\ge 1$. Noting that the sequence
$(X_i)_{i\ge
1}$ tends to infinity, it follows that
$$
\hat{w}(A,\utheta) \le {\delta +1\over w-\delta}.
$$
Choosing $\delta$ and $w$ arbitrarily close to  $0$ and to ${w}({}^tA)$
respectively, we obtain the second upper bound of (10).

\medskip

Furthermore, the preceding arguments enable us to establish Proposition
1.
The first assertion
follows immediately from the Theorem, since we just have to pick an
$n$-tuple
$\utheta$ out of a set of full Lebesgue measure. The proof of the
second assertion needs more work.
We begin by extracting some subsequence from the sequence of best
approximations $(\uy_i)_{i \ge 1}$.

We claim that there exists an increasing function
$\varphi : \bZ_{\ge 1} \to \bZ_{\ge 1}$ satisfying $\varphi (1) = 1$
and, for any integer $i \ge 2$,
$$
Y_{\varphi(i)} \ge (9n)^{1/2} Y_{\varphi (i-1)} \quad
{\rm and} \quad Y_{\varphi(i-1) + 1} \ge (9n)^{-1} Y_{\varphi (i)}.
\leqno (13)
$$
The function $\varphi$ is constructed in the following way.
Let $j > j'$ be two  indices such that
$Y_j \ge (9n)^{1/2} Y_{j-1}$ and $Y_{j'} \ge (9n)^{1/2} Y_{j'-1}$.
Suppose that $j'-1 = \varphi (h')$, and that the function $\varphi$
has already been defined for $1\le i \le  h'$.
We set $j-1 = \varphi (h)$ for some $h > h'$, which will be
specified later. We let
$\varphi(h-1)$ be the largest index $t \ge j'$ for which
$Y_{\varphi(h)} \ge (9n)^{1/2} \, Y_t$. We let
$\varphi(h-2)$ be the largest index $t \ge j'$ for which
$Y_{\varphi(h-1)} \ge (9n)^{1/2} \, Y_t$, and so on until it does not
exist any index $t$ as above.
We have just defined  $\varphi(h), \varphi(h-1), \ldots,
\varphi(h-h_0)$.
Then, we set $h = h_0 + h' + 1$, and we check that the
inequalities (13) are satisfied for $i= h'+1, \ldots , h_0 + h' + 1$.

This process does not apply when there are only finitely many
indices $j$ such that  $Y_j \ge (9n)^{1/2} Y_{j-1}$.
In this case, we denote by $g$ the largest of these indices ($g=1$ if
there
is none) and
we apply the above process to construct the initial values of
the function $\varphi$ up to  $g=\varphi(h)$.
   Next, we define
     $\varphi(h+1)$ as the smallest
    index $t$ for which $Y_t \ge (9n)^{1/2} \, Y_{\varphi(h)}$.
Then, we observe that  $Y_{\varphi(h+1)-1} < (9n)^{1/2} \,
Y_{\varphi(h)}$
and
$$
Y_{\varphi(h) + 1} \ge Y_{\varphi(h)} >
(9n)^{-1/2} \, Y_{\varphi(h+1)- 1} > (9n)^{-1} Y_{\varphi(h+1)},
$$
as required. We continue in this way, defining $\varphi(h+2)$ as the
smallest  index $t$ for which $Y_t \ge (9n)^{1/2} \, Y_{\varphi(h+1)}$,
and so on. The  inequalities (13) are then satisfied.

The first inequalities in (13) enable us to satisfy the assumptions
of  Lemma 2, page 86, from  [\Cas] for the sequence of
integer vectors $(\uy_{\varphi (i)})_{i \ge 1}$ with  $k=3$.
Consequently, there exists a real  $n$-tuple $\utheta$ such that
$$
\Vert y_{\varphi(i), 1} \theta_1 + \ldots + y_{\varphi(i), n} \theta_n
\Vert
\ge {1 \over 4}, \quad \hbox{for all $i \ge 1$}.  \leqno (14)
$$
Let  $\ux$ be a non-zero  integer $m$-tuple and let  $k$ be the index
defined
by the inequalities
$$
Y_{\varphi(k)} \le 9n (8m)^{m/n}  \, \vert \ux \vert^{m/n} <
Y_{\varphi(k+1)}.
$$
Taking into account (3), (5), (14) and (11) applied with
     $\uy = \uy_{\varphi (k)}$, we have
$$
{1 \over 4} \le (9 n^2) (8m)^{m/n} \, \vert \ux \vert^{m/n}
\, \Vert A\ux +\utheta\Vert +
m \vert \ux \vert \, Y_{\varphi(k)+1}^{-n/m}.
$$
By construction of the subsequence $(Y_{\varphi(i)})_{i \ge 1}$, we have
$Y_{\varphi(k)+1}^{-1} \cdot Y_{\varphi(k+1)} \le 9n$, so that
$$
{1 \over 4} \le (9 n^2) (8m)^{m/n} \, \vert \ux \vert^{m/n}
\, \Vert A\ux +\utheta\Vert +
m \, \bigl( 8m (9n)^{n/m} \bigr)^{-1} \, (9n)^{n/m},
$$
and
$$
\Vert A\ux +\utheta\Vert \ge {1 \over 72 n^2 (8m)^{m/n}} \,
\vert \ux \vert^{-m/n},
$$
as announced. \cqfd

\section
{5. The Corollary.}

Our Theorem reduces the determination of the measure of generic density
$$
\displaystyle \inf_{\utheta\in  \bR^n}w(A, \utheta) ={1\over
\hat{w}({}^tA)}
$$
of the group $\Gamma$ to the computation of the exponent
$\hat{w}({}^tA)$.
Any upper bound of $\hat{w}({}^tA)$
    implies a uniform lower bound for the exponents of approximation
$w(A,\utheta)$.
When
$$
A =(\xi, \dots , \xi^n),
$$
for a real transcendental number $\xi$, there are known upper bounds for
    $\hat{w}(A)$ and $\hat{w}({}^tA)$, which depend only upon  $n$.
Coming back to the specific notations of  [\BuLa]:
$$
\hat{w}(A) = \hat w_n(\xi) \and \hat{w}({}^tA) = \hat\lambda_n(\xi),
$$
we have
$$
\hat\lambda_n(\xi) \le {1\over \lceil n/2\rceil} \and
\hat w_n(\xi) \le 2n-1.
$$
The first upper bound is the main result of [\Lau], while
Theorem 2b of [\DaSc] is equivalent to the second one.
Actually, the results of [\Lau] and [\DaSc] are slightly sharper:
the refinement stated below Proposition 1 hold in these cases. Combined
with Lemma 3, they yield
the assertions (i) and (ii) of the  Corollary.

In  degree $n=2$, the exact upper bounds for the functions
$\hat\lambda_2(\xi)$ and
$\hat w_2(\xi)$ are known: Roy [\RoyA] and Arbour \& Roy
[\ArRo] have proved that
$$
\hat\lambda_2(\xi) \le {\sqrt{5}-1\over 2} \and
\hat w_2(\xi) \le {\sqrt{5} +3\over 2},
$$
and both equalities hold  when $\xi$ is a Fibonacci continued fraction.
The assertion (iii) of the Corollary is then the translation of our
Theorem in this particular case.

\medskip
\noindent
{\bf Remarks.}
(i)
The exponents $\hat\lambda_2(\xi)$ and $\hat w_2(\xi)$ have been
computed more generally in  [\BuLa]  for any Sturmian continued fraction
     $\xi$ of irrational angle  $\varphi$. It turns out that
$$
\hat\lambda_2(\xi) >{1\over 2}  \and  \hat w_2(\xi) > 2
$$
when the partial quotients  in the continued fraction expansion
of  $\varphi$ are bounded.
The associated subgroups $\Gamma$
$$
\bZ + \bZ \xi +\bZ \xi^2 \and \bZ\pmatrix{\xi\cr\xi^2\cr} +\bZ^2
$$
are then  dense in  $\bR$ and $\bR^2$ respectively, and their generic
exponents
of density are less than $2$ and $1/2$, respectively.

\medskip

(ii) When  $\xi$ is some real number  connected with the
      Fibonacci continued fraction, Roy [\RoyB] has proved that we have
a lower bound
of the form
$$
     \vert P(\xi)+\xi^3\vert \gg H(P)^{-(1+\sqrt{5})/2}
$$
for any quadratic polynomial   $P(X)$ with integer coefficients.
This means that the number $\theta = \xi^3$
shares the almost sure property stated in
the part (iii) of the Corollary. The difficult point in Roy's proof
is to  verify that a lower bound similar to (14) is valid for
$\utheta = (\xi^3)$.

\section
{6.  Hausdorff dimension.}

The   Hausdorff dimension is a useful tool to
discriminate between sets of
Lebesgue measure zero and thus to prove the existence of real numbers
having fine properties of Diophantine approximation
    (see, for instance, Chapter 5 of [\BuLiv]). We denote it by $\dim$ 
and
direct  the reader   to the books of Falconer [\FalA, \FalB]
for the definition and the properties of Hausdorff dimension.

Let  $A$ be a matrix  in  $M_{n, m} (\bR)$ and let
$w > 1/\hat{w}({}^tA)$ be a real number.
In  view of our  Theorem, it is  natural to ask whether there
exists a real  $n$-tuple $\utheta$ such that
$w(A, \utheta) = w$. As a  first step, we wish to  compute the
    Hausdorff dimension of the null set
$$
\cU_w (A)=\biggl\{\utheta \in \bR^{n} : \Vert A\ux + \utheta \Vert \leq
{ 1 \over |\ux|^{w}} \hbox{ for infinitely many }\ux
\hbox{ in }\bZ^{m}\biggr\}.
$$
This  question has been solved  in the simplest case when
    $A = (\xi)$, for any given irrational number  $\xi$, by
Bugeaud [\BuGl] and, independently, by Schmeling \&
Troubetzkoy [\ScTr]. We then have $w((\xi), \utheta) = 1$ for almost all
    real numbers $\theta$ and the
Hausdorff dimension of the set $\cU_w ((\xi))$ is  equal to  $1/w$,
for any $w\ge 1$.

This question remains unanswered for other
    matrices $A$. Nevertheless, Bugeaud \& Chevallier [\BuCh]
have proved that,
for almost all matrices $A$ in  $M_{n, m} (\bR)$
and for any real number $w \ge m/n$, we have
$$
\dim \, \cU_w (A) = \dim \, \{\utheta \in \bR^{n} : w(A, \utheta) \ge w
\}
= {m \over w}.
$$
Besides,  Theorem 3 of [\BuCh] asserts that, if  $A$ is a  column
matrix, then
$$
\dim \, \cU_w (A) = \dim \, \{\utheta \in \bR^{n} : w(A, \utheta) \ge w
\}
= {1 \over w}
$$
for any real number $w \ge 1$. In the above examples, the
sets $\{\utheta \in \bR^{n} : w(A, \utheta) \ge w \}$
and $\cU_w (A)$ have
the same  Hausdorff dimension (although the first one contains the
second).

The results of
[\BuCh] indicate that the situation is much more complicated when the
matrix $A$ is not of the form  $(\xi)$.
It seems to us that the determination  in our general framework
of the  Hausdorff dimension of the sets $\cU_w (A)$
is a quite difficult problem.
Nonetheless, it is possible to show that this dimension is strictly
less than
    $n$, whenever $w > 1/\hat{w}({}^tA)$.

\proclaim Proposition 2. Let  $A$ be a  matrix in $M_{n, m} (\bR)$
and let  $w$  be a real number $> 1/\hat{w}({}^tA)$. The set
$$
\{\utheta \in \bR^{n} : w(A, \utheta) \ge w \}
$$
has Lebesgue measure zero, and its  Hausdorff dimension is
    strictly less than ~$n$.

\proof Set
$$
\delta = {w \hat{w}({}^tA) - 1 \over
2 + w + 1/ \hat{w}({}^tA)}
$$
and notice that the inequality  (6) determines
in the hypercube $[0, 1]^n$ a subset contained in the union of at most
    $c(n) \, Y_i \cdot Y_i^{(\delta + 1) (n - 1)}$ hypercubes
with edge $Y_i^{-\delta - 1}$, where $c(n)$ denotes some suitable
constant,  depending only upon  $n$. Since, by Lemma 1, the series
$$
\sum_{i \ge 1} \, Y_i^{1 + (\delta + 1) (n - 1)} \cdot Y_i^{-(\delta +
1)s}
$$
converges for any  $s > n-1 + 1/(\delta + 1)$, the
Hausdorff--Cantelli Lemma (cf. for example [\BuLiv], Chapter 5) ensures
us
that the  Hausdorff dimension of the set
$$
\cV_\delta =
\{ \utheta \in \bR^{n} :
\Vert y_{i,1}\theta_1 + \cdots + y_{i,n}\theta_n\Vert  <
Y_i^{-\delta} \quad \hbox{for infinitely many $i$} \}
$$
is  bounded from above by
$ n-1 + 1/(\delta + 1)$. This is strictly less than $n$
since  $\delta$ is positive.
Let $\utheta$ be in the complement of $\cV_{\delta}$, and follow
again the proof of our Theorem. The inequality (12)
is then  satisfied for any sufficiently large integer $i$. Thus, we have
the upper bound
$$
w(A, \utheta) \le {\delta + 1 \over \hat{w}({}^tA) - \delta},
$$
and, by our choice of  $\delta$,
$$
w(A, \utheta) \le {1 \over 2} \biggl( {1 \over \hat{w}({}^tA)} + w
\biggr)
= w - {1 \over 2} \biggl( w - {1 \over \hat{w}({}^tA)} \biggr).
$$
Consequently, the set $\cU_w ( A)$ is contained in $\cV_\delta $.
This remark  concludes the proof of Proposition 2. \cqfd

In view of the results of [\RoyA, \RoyB, \BuLa], it may be possible
to determine the Hausdorff dimensions of the sets $\cU_w (A)$
and $\cU_w ({}^t A)$, when  $A = (\xi, \xi^2)$ and
$\xi$ is a Sturmian continued fraction. We plan to return to these
questions
later. Let us simply remark that our
Proposition 2 implies that  $\dim \, \cU_2 (A) < 1$ and
$\dim \, \cU_{1/2} ({}^t A) < 2$, unlike in the generic
situation.

\vskip 8mm

\centerline{\bf References }

\vskip 5mm

\item{[\ArRo]}
B. Arbour and D. Roy,
{\it A Gel'fond type criterion in degree two},
Acta Arith. 111 (2004), 97--103.

\item{[\BuGl]}
Y. Bugeaud,
{\it A note on inhomogeneous Diophantine approximation},
Glasgow Math. J. 45 (2003), 105--110.

\item{[\BuLiv]}
Y. Bugeaud,
Approximation by algebraic numbers,
Cambridge Tracts in Mathematics, Cambridge, 2004.

\item{[\BuCh]}
Y. Bugeaud and N. Chevallier,
{\it On simultaneous inhomogeneous Diophantine approximation},
Preprint.

\item{[\BuLa]}
Y. Bugeaud and M. Laurent,
{\it Exponents of Diophantine Approximation and
Sturmian Continued Fractions}, Preprint.

\item{[\Cas]}
      J. W. S. Cassels,
An introduction to Diophantine Approximation,
Cambridge Tracts in Math. and Math. Phys., vol. 99, Cambridge
University Press, 1957.

\item{[\CasB]}
      J. W. S. Cassels,
An introduction to the Geometry of Numbers,
Springer Verlag, 1997.

\item{[\DaSc]}
       {H. Davenport and W. M. Schmidt},
{\it Approximation to real numbers by
algebraic integers}, Acta Arith. {15} (1969), 393--416.

\item{[\DSa]}
H. Davenport and W. M. Schmidt,
{\it Dirichlet's theorem on
Diophantine approximation}, Symposia Mathematica, Vol. IV
(INDAM, Rome, 1968/69), pp. 113--132, Academic Press, London, 1970.

\item{[\FalA]}
K. Falconer,
The geometry of fractal sets,
Cambridge Tracts in Mathematics 85, Cambridge University Press, 1985.

\item{[\FalB]}
K. Falconer,
Fractal Geometry : Mathematical Foundations and
Applications, John Wiley \& Sons, 1990.

\item{[\Gro]}
A. V. Groshev,
{\it A theorem on systems of linear forms},
Dokl. Akad. Nauk SSSR 19 (1938), 151--152 (in Russian).

\item{[\Jar]}
V. Jarn\'\i k,
{\it Eine Bemerkung \"uber diophantische Approximationen},
Math. Z. 72 (1959), 187--191.

\item{[\JarB]}
V. Jarn\'\i k,
{\it Contributions to the theory  of homogeneous linear  diophantine 
approximations},
Czech. Math. J.4 (79) (1954), 330--353 (in Russian, French summary).

\item{[\Khi]}
A. Ya. Khintchine,
{\it \"Uber eine Klasse linearer diophantischer Approximationen},
Rendiconti Circ. Mat. Palermo 50 (1926), 170--195.

\item{[\LagA]}
       {J. C. Lagarias},
{\it Best simultaneous diophantine approximations. I.
Growth rates of best approximation denominators},
Trans. Amer. Math. Soc.  {272} (1982), 545--554.

\item{[\LagB]}
       {J. C. Lagarias},
{\it Best Diophantine approximations to a set of linear forms},
J. Austral. Math. Soc. Ser. A  {34} (1983), 114--122.

\item{[\Lau]}
       M. Laurent,
{\it Simultaneous rational approximation to the successive powers
of a real number},
Indag. Math. 11 (2003), 45--53.

\item{[\RoyA]}
       D. Roy,
{\it Approximation to real numbers by cubic algebraic numbers, I},
Proc. London Math. Soc. 88 (2004), 42--62.

\item{[\RoyB]}
       D. Roy,
{\it Approximation to real numbers by cubic algebraic numbers, II},
Annals of Math. 158 (2003), 1081--1087.

\item{[\ScTr]}
J. Schmeling and S. Troubetzkoy,
{\it Inhomogeneous Diophantine Approximation and Angular Recurrence for
Polygonal Billiards},
Mat. Sbornik 194 (2003), 295--309.

\item{[\Schm]}
W. M. Schmidt,
Diophantine Approximation, Lecture Notes in Math.
{785}, Springer, Berlin, 1980.

\item{[\Wa]}
{M. Waldschmidt},
Topologie des points rationnels. Cours de D.E.A. Paris.

\vskip1cm

\noindent Yann Bugeaud  \hfill{Michel Laurent}

\noindent Universit\'e Louis Pasteur
\hfill{Institut de Math\'ematiques de Luminy}

\noindent U. F. R. de math\'ematiques
\hfill{C.N.R.S. -  U.P.R. 9016 - case 907}

\noindent 7, rue Ren\'e Descartes      \hfill{163, avenue de Luminy}

\noindent 67084 STRASBOURG  (FRANCE)
\hfill{13288 MARSEILLE CEDEX 9  (FRANCE)}

\vskip2mm

\noindent {\tt bugeaud@math.u-strasbg.fr}
\hfill{{\tt laurent@iml.univ-mrs.fr}}

\end